\documentclass[titlepage,draft,12pt]{article} 
\usepackage{amsfonts,latexsym,pstricks,pst-plot} 
\usepackage[a4paper]{geometry}
\newcommand{\GM}{\mathcal{GM}}
\newcommand{\ocont}{$\omega$-continuous}
\newcommand{\G}{\mathcal{G}}
\newcommand{\trebar}[1]{\|#1\|}
\newcommand{\lk}{\lambda_{k}}
\newcommand{\ep}{\varepsilon}
\newcommand{\re}{{\mathbb{R}}}
\newcommand{\n}{{\mathbb{N}}}
\newcommand{\m}[1]{m(\au{#1}^{2})}
\newcommand{\au}[1]{|A^{1/2}#1|}

\newtheorem{thm}{Theorem}[section]

\newtheorem{prop}[thm]{Proposition}
\newtheorem{defn}[thm]{Definition}

\newtheorem{open}{Open problem}

\title{Kirchhoff equations in generalized Gevrey spaces: local
existence, global existence, uniqueness}
\author{Marina Ghisi\medskip\\ 
{\normalsize Universit\`a degli Studi di Pisa} \\
{\normalsize Dipartimento di Matematica ``Leonida Tonelli''}\\
{\normalsize PISA (Italy)}\\
{\normalsize e-mail: \texttt{ghisi@dm.unipi.it}}
\and 
Massimo Gobbino\medskip\\ 
{\normalsize Universit\`a degli Studi di Pisa} \\
{\normalsize Dipartimento di Matematica Applicata ``Ulisse Dini''}\\ 
{\normalsize e-mail: \texttt{m.gobbino@dma.unipi.it}}}
\date{}

\begin{document}
\maketitle
\begin{abstract}
	In this note we present some recent results for Kirchhoff
	equations in generalized Gevrey spaces.  We show that these spaces
	are the natural framework where classical results can be unified
	and extended.  In particular we focus on existence and uniqueness
	results for initial data whose regularity depends on the
	continuity modulus of the nonlinear term, both in the strictly
	hyperbolic case, and in the degenerate hyperbolic case.
	
\vspace{1cm}

\noindent{\bf Mathematics Subject Classification 2000 (MSC2000):}
35L70, 35L80, 35L90.

\vspace{1cm} 

\end{abstract}
 
\section{Introduction}

Let $H$ be a separable real Hilbert space.  For every $x$ and $y$ in
$H$, $|x|$ denotes the norm of $x$, and $\langle x,y\rangle$ denotes
the scalar product of $x$ and $y$.  Let $A$ be a self-adjoint linear
operator on $H$ with dense domain $D(A)$.  We assume that $A$ is
nonnegative, namely $\langle Ax,x\rangle\geq 0$ for every $x\in D(A)$,
so that for every $\alpha\geq 0$ the power $A^{\alpha}x$ is defined
provided that $x$ lies in a suitable domain $D(A^{\alpha})$.

Given a continuous function $m:[0,+\infty)\to[0,+\infty)$ we consider 
the Cauchy problem
\begin{equation}
	u''(t)+\m{u(t)}Au(t)=0 
	\hspace{2em}\forall t\in[0,T),
	\label{pbm:h-eq}
\end{equation}
\begin{equation}
	u(0)=u_0,\hspace{3em}u'(0)=u_1.
	\label{pbm:h-data}
\end{equation}
 
It is well known that (\ref{pbm:h-eq}), (\ref{pbm:h-data}) is the
abstract setting of the Cauchy-boundary value problem for the
quasilinear hyperbolic integro-differential partial differential
equation
\begin{equation}
	u_{tt}(t,x)-
	m{\left(\int_{\Omega}\left|\nabla u(t,x)\right|^2\,dx\right)}
	\Delta u(t,x)=0
	\hspace{2em}
	\forall(x,t)\in\Omega\times[0,T),
	\label{eq:k}
\end{equation}
where $\Omega\subseteq\re^{n}$ is an open set, and $\nabla u$ and
$\Delta u$ denote the gradient and the Laplacian of $u$ with respect
to the space variables. Let us set
$$\mu:=\inf_{\sigma\geq 0}m(\sigma).$$

Equation (\ref{pbm:h-eq}) or (\ref{eq:k}) are called \emph{strictly
hyperbolic} if $\mu>0$, and \emph{weakly} (or \emph{degenerate})
\emph{hyperbolic} if $\mu=0$.

From the mathematical point of view, (\ref{eq:k}) is probably the
simplest example of quasilinear hyperbolic equation.  From the
mechanical point of view, this Cauchy boundary value problem is a
model for the small transversal vibrations of an elastic string
($n=1$) or membrane ($n=2$).  In this context it was introduced by
G.\ Kirchhoff in \cite{kirchhoff}.

We refer to \cite{A-trieste} for a sketch of the deduction of
(\ref{eq:k}) from the system of (local) equations of elasticity, and 
to \cite{ap} for the standard arguments in functional analysis 
leading from (\ref{eq:k}) to (\ref{pbm:h-eq}).

This equation has generated a considerable literature after the
pioneering paper by S.\ Bernstein~\cite{bernstein}.  For more
details on previous works we refer to the introductions of the
following sections.  In this note we focus on the basic analytic
questions for a partial differential equation, namely local existence,
uniqueness, and global existence.

Section~\ref{sec:local} is devoted to local existence.  We show that a
local solution of (\ref{pbm:h-eq}) exists provided that the initial
data (\ref{pbm:h-data}) are regular enough, depending on the
continuity modulus of $m$.  This result is an interpolation between
the two extreme cases, namely the classical local existence results
for analytic data and continuous $m$, or for Sobolev data and
Lipschitz continuous $m$.  We show that these local solutions satisfy
the expected properties of propagation of regularity and continuous
dependence on the data.  More important, we show with some
counterexamples that the spaces involved in the local existence
results are optimal.

Section~\ref{sec:uniqueness} is devoted to uniqueness issues.  We
present a uniqueness result in which the nonlinear term is \emph{not}
required to be Lipschitz continuous.

Section~\ref{sec:global} is devoted to global existence.  After
reviewing the classical global existence results, we state our result
concerning ``spectral gap'' initial data.  This special class of
initial data is small in the sense that it is not a vector space, and
it does not even contain all analytic data, but it is large in the
sense that every initial condition in the spaces involved in the local
existence result is the sum of two initial data for which the solution
is actually global.  In particular, a solution can globally exist even
if its initial data have only the minimal regularity required by the
local existence result.

Finally, Section~\ref{sec:open} is devoted to open problems. We 
recall some old and new unsolved questions which should inspire the 
future investigations in this challenging research field.

We conclude by pointing out that there is a considerable literature on
equation (\ref{pbm:h-eq}) or (\ref{eq:k}) with additional
\emph{dissipative terms}.  The interested reader is referred to the
survey~\cite{gg:survey-diss}.

\setcounter{equation}{0}
\section{Local existence}\label{sec:local}

\paragraph{Previous works}

Local existence has been proved in the last century with two opposite 
sets of assumptions.
\begin{enumerate}
	\item[(A)] The case where $m$ is Lipschitz continuous, the
	equation is strictly hyperbolic, and initial data have Sobolev
	regularity.  Under these assumptions a local existence result was
	first proved by S.\ Bernstein~\cite{bernstein}, and then
	extended with increasing generality by many authors.  The more
	general form was probably stated by A.\ Arosio and
	S.\ Panizzi in~\cite{ap}, where they proved that problem
	(\ref{pbm:h-eq}), (\ref{pbm:h-data}) is well posed in the phase
	space $D(A^{3/4})\times D(A^{1/4})$.

	\item[(B)] The case where $m$ is continuous, the equation is
	weakly hyperbolic, and initial data are analytic.  A local (and
	actually global, see Section~\ref{sec:global}) existence result
	under these assumptions was proved with increasing generality by
	S.\ Bernstein~\cite{bernstein}, S. I. Pohozaev~\cite{poho-an}, A.\
	Arosio and S.\ Spagnolo~\cite{as}, P.\ D'Ancona and S.\
	Spa\-gno\-lo~\cite{das-an-1,das-an-2}.
\end{enumerate}

More recently, F.\ Hirosawa~\cite{hirosawa-main} considered
equation~(\ref{eq:k}) with $\Omega=\re^{n}$, and proved a local
existence result in classes of initial data depending on the
continuity modulus of $m$.  The rough idea is that the less regular is
the nonlinear term, the more regularity is required on initial data.
This result interpolates (A) and (B).

Our contribution is twofold: in \cite{gg:local} we extended these
intermediate results from the special concrete case $\Omega=\re^{n}$
to the general abstract setting, and we provided counterexamples in
order to show their optimality.  Let us introduce the functional
setting needed in these statements.

\paragraph{Functional spaces and continuity moduli}

For the sake of simplicity we assume that $H$ admits a countable
complete orthonormal system $\{e_{k}\}_{k\geq 1}$ made by eigenvectors
of $A$.  We denote the corresponding eigenvalues by $\lambda_{k}^{2}$
(with $\lk\geq 0$), so that $Ae_{k}=\lambda_{k}^{2}e_{k}$ for every
$k\geq 1$.  Every $u\in H$ can be written in a unique way in the form
$u=\sum_{k=1}^{\infty}u_{k}e_{k}$, where $u_{k}=\langle
u,e_{k}\rangle$ are the components of $u$.  In other words, every
$u\in H$ can be identified with the sequence $\{u_{k}\}$ of its
components, and under this identification the operator $A$ acts
component-wise by multiplication.

This simplifying assumption is not so restrictive.  Indeed the
spectral theorem for self-adjoint unbounded operators on a separable
Hilbert space (see \cite[Chapter~VIII]{reed}) states that any such
operator is unitary equivalent to a multiplication operator on some
$L^{2}$ space.
		
More precisely, for every $H$ and $A$ there exist a measure space
$(M,\mu)$, a function $a(\xi)\in L^{2}(M,\mu)$, and a unitary operator
$H\to L^{2}(M,\mu)$ which associates to every $u\in H$ a function
$f(\xi)\in L^{2}(M,\mu)$ in such a way that $Au$ corresponds to the
product $a(\xi)f(\xi)$.
		
As a consequence, all the spaces we define in terms of $u_{k}$ and
$\lambda_{k}$ can be defined in the general case by replacing the
sequence of components $\{u_{k}\}$ of $u$ with the function $f(\xi)$
corresponding to $u$, the sequence $\{\lambda_{k}\}$ of eigenvalues of
$A$ with the function $a(\xi)$, and summations over $k$ with integrals
over $M$ in the variable $\xi$ with respect to the measure $\mu$.
Similarly, there is no loss of generality in proving our existence and
uniqueness results (Theorems~\ref{thm:local}, \ref{thm:uniqueness},
\ref{thm:global}) using components, as we did
in~\cite{gg:local,gg:uniqueness,gg:global}.  On the contrary,
existence of countably many \emph{eigenvalues} is essential in the
construction of our counterexamples, as stated in
Theorem~\ref{thm:counterex-ndg} and Theorem~\ref{thm:counterex-dg}.

Coming back to functional spaces, using components we have that 
$$D(A^{\alpha}):=\left\{u\in H:\sum_{k=1}^{\infty}
\lambda_{k}^{4\alpha}u_{k}^{2}<+\infty\right\},
\hspace{3em}
D(A^{\infty}):=\bigcap_{\alpha\geq 0}^{}D(A^{\alpha}).$$

Let now $\varphi:[0,+\infty)\to[1,+\infty)$ be any function.  Then for
every $\alpha\geq 0$ and $r>0$ one can set
\begin{equation}
	\trebar{u}_{\varphi,r,\alpha}^{2}:=\sum_{k=1}^{\infty}\lambda_{k}^{4\alpha}
	u_{k}^{2} \exp\left(\strut r\varphi(\lambda_{k})\right),
	\label{defn:trebar}
\end{equation}
and then define the \emph{generalized Gevrey spaces} as
\begin{equation}
	\G_{\varphi,r,\alpha}(A):=
	\left\{u\in H:\trebar{u}_{\varphi,r,\alpha}<+\infty\right\}.
	\label{defn:G}
\end{equation}

These spaces can also be seen as the domain of the operator
$A^{\alpha}\exp\left(r\varphi(A)\right)$, and they are a
generalization of the usual spaces of Sobolev, Gevrey or analytic
functions.  They are Hilbert spaces with norm
$(|u|^{2}+\trebar{u}_{\varphi,r,\alpha}^{2})^{1/2}$, and they form a
scale of Hilbert spaces with respect to the parameter $r$.

A \emph{continuity modulus} is a continuous increasing function
$\omega:[0,+\infty)\to[0,+\infty)$ such that $\omega(0)=0$, and
$\omega(a+b)\leq\omega(a)+\omega(b)$ for every $a\geq 0$ and $b\geq
0$.

The function $m$ is said to be $\omega$-continuous if there exists a
constant $L\in\re$ such that
\begin{equation}
	|m(a)-m(b)|\leq
	L\,\omega(|a-b|)
	\hspace{3em}
	\forall a\geq 0,\ \forall b\geq 0.
	\label{hp:ocont}
\end{equation}

We point out that the set of $\omega$-continuous functions depends
only on the values of $\omega(\sigma)$ for $\sigma$ in a neighborhood
of $0$, while when $\varphi$ is continuous the space
$\G_{\varphi,r,\alpha}(A)$ depends only on the values of
$\varphi(\sigma)$ for large $\sigma$.

\paragraph{Our local existence results}

The following statement sums up the state of the art concerning
existence of local solutions.

\begin{thm}[Local existence]\label{thm:local}
	Let $H$ be a separable Hilbert space, and let $A$ be a nonnegative
	self-adjoint (unbounded) operator on $H$ with dense domain.  Let
	$\omega$ be a continuity modulus, let
	$m:[0,+\infty)\to[0,+\infty)$ be an \ocont\ function, and let
	$\varphi:[0,+\infty)\to[1,+\infty)$.
	
	Let us assume that there exists a constant $\Lambda$ such that 
	either
	\begin{equation}
		\sigma
		\omega\left(\frac{1}{\sigma}\right)\leq\Lambda\varphi(\sigma)
		\quad\quad
		\forall\sigma> 0
		\label{hp:phi-ndg}
	\end{equation}
	in the strictly hyperbolic case, or 
	\begin{equation}
		\sigma\leq\Lambda\varphi\left(\frac{\sigma}{
		\sqrt{\omega(1/\sigma)}}\right)
		\quad\quad
		\forall\sigma> 0
		\label{hp:phi-dg}
	\end{equation}
	in the weakly hyperbolic case.
	Let us finally assume that
	\begin{equation}
		(u_{0},u_{1})\in\G_{\varphi,r_{0},\alpha+1/2}(A)\times
		\G_{\varphi,r_{0},\alpha}(A)
		\label{hp:data}
	\end{equation}
	for some $r_{0}>0$, and some $\alpha\geq 1/4$.
	
	Then there exist $T>0$, and $R>0$ with $RT<r_{0}$ such that
	problem (\ref{pbm:h-eq}), (\ref{pbm:h-data}) admits at least one
	local solution $u(t)$ in the space
	\begin{equation}
			 C^{1}\left([0,T];\G_{\varphi,r_{0}-Rt,\alpha}(A)\right)\cap
			C^{0}\left([0,T];\G_{\varphi,r_{0}-Rt,\alpha+1/2}(A)\right).
		\label{th:reg-sol}
	\end{equation}
\end{thm}

Condition (\ref{th:reg-sol}), with the range space depending on
time, simply means that 
$$u\in
C^{1}\left([0,\tau];\G_{\varphi,r_{0}-R\tau,\alpha}(A)\right)\cap
C^{0}\left([0,\tau];\G_{\varphi,r_{0}-R\tau,\alpha+1/2}(A)\right)$$
for every $\tau\in(0,T]$.  This amounts to say that scales of Hilbert
spaces are the natural setting for this problem.

Admittedly assumptions (\ref{hp:phi-ndg}) and (\ref{hp:phi-dg}) do not
lend themselves to a simple interpretation.  The basic idea is that in
the strictly hyperbolic case the best choice for $\varphi$, namely the
choice giving the largest space of initial data, is always
$\varphi(\sigma)=\sigma\omega(1/\sigma)$.  In the weakly hyperbolic
case things are more complex because condition (\ref{hp:phi-dg}) is
stated in an implicit form.  In this case the best choice for
$\varphi$ is the inverse of the function
$\sigma\to\sigma/\sqrt{\omega(1/\sigma)}$.  Note that this inverse
function is always $o(\sigma)$ as $\sigma\to +\infty$.
Tables~\ref{tab:local-sh} and~\ref{tab:local-wh} provide examples of
pairs of functions $m$, $\varphi$ satisfying (\ref{hp:phi-ndg}) and
(\ref{hp:phi-dg}).

\begin{table}[htb]
	\centering
	\renewcommand{\arraystretch}{1.4}
	\begin{tabular}{|c|c||c|c|}
		\hline
		$\omega(\sigma)=$ & $m$ is \ldots & $\varphi(\sigma)=$ & Local
		existence with data in \ldots \\
		\hline\hline
		any & continuous & $\sigma$ & analytic functions (never
		optimal) \\
		\hline
		any & continuous & $\sigma\omega(1/\sigma)$ & space larger
		than analytic functions \\
		\hline
		$\sigma^{\beta}$ & $\beta$-H\"{o}lder cont. &
		$\sigma^{1-\beta}$ & Gevrey space of order $(1-\beta)^{-1}$ \\
		\hline
		$\sigma|\log\sigma|$ & Log-Lipschitz cont. & $\log\sigma$ & 
		$D(A^{\alpha+1/2})\times D(A^{\alpha})$ with $\alpha>1/4$  \\
		\hline
		$\sigma$ & Lipschitz cont. & 1 & $D(A^{3/4})\times D(A^{1/4})$  \\
		\hline
	\end{tabular}
	\caption{Examples of relations between the regularity of $m$ and
	the regularity of initial data for local existence in the
	\emph{strictly} hyperbolic case}
	\label{tab:local-sh}
\end{table}

\begin{table}[htb]
	\centering
	\renewcommand{\arraystretch}{1.4}
	\begin{tabular}{|c|c||c|c|}
		\hline
		$\omega(\sigma)=$ & regularity of $m$ & $\varphi(\sigma)=$ &
		Local existence with data in \ldots \\
		\hline\hline
		any & continuous & $\sigma$  & analytic functions (never optimal)  \\
		\hline
		any & continuous & $o(\sigma)$  & space larger than analytic functions  \\
		\hline
		$\sigma^{\beta}$ & $\beta$-H\"{o}lder cont.  &
		$\sigma^{2/(\beta+2)}$ & Gevrey space of order $1+\beta/2$ \\
		\hline
		$\sigma$ & Lipschitz cont.  & $\sigma^{2/3}$ & Gevrey space of
		order $3/2$ \\
		\hline
	\end{tabular}
	\caption{Examples of relations between the regularity of $m$ and
	the regularity of initial data for local existence in the
	\emph{weakly} hyperbolic case}
	\label{tab:local-wh}
\end{table}

As one could easily expect, assumption (\ref{hp:phi-dg}) is always
stronger than assumption (\ref{hp:phi-ndg}).  We remark that, since we
are interested in local solutions, assumption (\ref{hp:phi-ndg}) is
the relevant one also when the equation is degenerate but the initial
condition $u_{0}$ satisfies the nondegeneracy condition
$$m(|A^{1/2}u_{0}|^{2})\neq 0.$$

In this case problem (\ref{pbm:h-eq}), (\ref{pbm:h-data}) is called
\emph{mildly degenerate}.  As observed in~\cite{ag}, in this situation
it is enough to solve the problem with a different nonlinearity which
is strictly hyperbolic and coincides with the given $m$ in a
neighborhood of $|A^{1/2}u_{0}|^{2}$.  The solution of the modified
problem is thus a solution of the original problem for $t$ small
enough.

Assumption (\ref{hp:phi-dg}) is therefore the relevant one only when
$m(|A^{1/2}u_{0}|^{2})=0$.  This is usually called the \emph{really
degenerate} case.

The proof of Theorem~\ref{thm:local} relies on standard techniques. 
The first step is remarking that (\ref{pbm:h-eq}) admits a 
first-order conserved energy, namely the Hamiltonian
\begin{equation}
	\mathcal{H}(t):=|u'(t)|^{2}+M(|A^{1/2}u(t)|^{2}),
	\label{defn:H}
\end{equation}
where $M(\sigma)$ is any function such that $M'(\sigma)=m(\sigma)$ 
for every $\sigma\geq 0$. This is the reason why $D(A^{1/2})\times H$ 
is called the \emph{energy space}.

The second step is to consider the linearization of (\ref{pbm:h-eq}), 
namely equation
\begin{equation}
	u''(t)+c(t)Au(t)=0,
	\label{eq:linear}
\end{equation}
where now $c(t):=m(|A^{1/2}u(t)|^{2})$ is thought as a given
coefficient.  The theory of such linear hyperbolic equations with
time-dependent coefficients was developed by F.\ Colombini, E.\ De
Giorgi and S.\ Spagnolo~\cite{dgcs} in the strictly hyperbolic case,
and by F.\ Colombini, E.\ Jannelli and S.\ Spagnolo~\cite{cjs} in the
weakly hyperbolic case. The result is that a local solution exists 
provided that the regularity of the initial data is related to the 
continuity modulus of $c(t)$ as in Theorem~\ref{thm:local}. 

Unfortunately the boundedness of the Hamiltonian (\ref{defn:H}) is 
not enough to control the oscillations of $c(t)$. The main point is 
therefore to prove an a priori estimate for
$$\frac{d}{dt}|A^{1/2}u(t)|^{2}= 2\langle
A^{1/2}u(t),A^{1/2}u'(t)\rangle= 2\langle
A^{3/4}u(t),A^{1/4}u'(t)\rangle,$$
which in turn is achieved through an a priori estimate for the higher 
order energy
$$|A^{1/4}u'(t)|^{2}+|A^{3/4}u(t)|^{2}.$$

This a priori estimate provides an a priori control on the continuity
modulus of $c(t)$.  One can therefore apply the linear theory and
obtain, for example, the so called \emph{propagation of regularity},
namely the fact that solutions belong to the same space (or more
precisely to the same scale of spaces) of the initial condition.  The
precise statement is the following.

\begin{thm}[A priori estimate, Propagation of
regularity]\label{thm:apriori-est}
	Let $H$, $A$, $\omega$, $m$, $\varphi$, $\Lambda$, $u_{0}$,
	$u_{1}$, $r_{0}$, $\alpha$ be as in Theorem~\ref{thm:local}.
	
	Then there exist positive real numbers $T$, $K$, $R$, with
	$RT<r_{0}$, such that every solution
	\begin{equation}
		u\in C^{1}\left([0,T];D(A^{1/4})\right)\cap
		C^{0}\left([0,T];D(A^{3/4})\right)	
		\label{hp:reg-u}
	\end{equation}
	of problem (\ref{pbm:h-eq}), (\ref{pbm:h-data}) satisfies
	$$|A^{1/4}u'(t)|^{2}+|A^{3/4}u(t)|^{2}\leq K 
	\quad\quad
	\forall t\in[0,T],$$
	and actually $u$ belongs to the space (\ref{th:reg-sol}).
\end{thm}

The constants $T$, $K$, $R$ depend only on $\omega$, $m$, and on the
norms of $u_{0}$ and $u_{1}$ in the spaces
$\G_{\varphi,r_{0},\alpha+1/2}(A)$ and $\G_{\varphi,r_{0},\alpha}(A)$,
respectively.

In Theorem~\ref{thm:apriori-est}, as in every a priori estimate, we
assumed the existence of a solution.  So the final step in the proof
of Theorem~\ref{thm:local} is proving that a solution exists.  Thanks
to the a priori estimate this can be done in several standard ways.

A first possibility is to use Galerkin approximations.  In this case
one approximates $(u_{0},u_{1})$ with a sequence of data
$(u_{0n},u_{1n})$ belonging to $A$-invariant subspaces of $H$ where
the restriction of $A$ is a \emph{bounded} operator.  For such data
solutions exist, and thanks to the a priori estimate the corresponding
coefficients $c_{n}(t)$ are relatively compact in $C^{0}([0,T])$.  The
conclusion follows from the fact that solutions of the linear problem
depend continuously on the initial data and on the coefficient $c(t)$ 
(see~\cite{dgcs} and~\cite{cjs}).

A second possibility is to apply Schauder's fixed point theorem in the
space of coefficients.  In this case one defines $X_{K,T}$ as the
space of functions $a:[0,T]\to\re$ such that
$a(0)=|A^{1/2}u_{0}|^{2}$, and with Lipschitz constant less than or
equal to $K$.  For every $a\in X_{K,T}$, one considers the solution
$u(t)$ of the linear problem (\ref{eq:linear}) with $c(t)=m(a(t))$,
and initial data (\ref{pbm:h-data}).  Finally one sets
$[\Phi(a)](t):=|A^{1/2}(t)|^{2}$.  For suitable values of $K$ and
$T$ (those given by Theorem~\ref{thm:apriori-est}), $\Phi$ defines a
continuous map from $X_{K,T}$ into itself which has a fixed point due
to Schauder's theorem.  This fixed point corresponds to a solution of
(\ref{pbm:h-eq}), (\ref{pbm:h-data}).

The same techniques (a priori estimate + compactness + results for 
the linear equation) lead to the continuity of the map
$$(\mbox{initial data}, m)\to\mbox{solution}.$$

Since the solution is not necessarily unique, this has to be intended
in the sense that ``the limit of solutions is again a solution''. The 
precise statement is the following.

\begin{thm}[Continuous dependence on initial data]\label{thm:cont-dep}
	Let $H$, $A$, $\omega$, $m$, $\varphi$, $\Lambda$, $u_{0}$,
	$u_{1}$, $r_{0}$, $\alpha$, $T$, $R$ be as in
	Theorem~\ref{thm:apriori-est}. 
	
	Let $\{m_{n}\}$ be a sequence of $\omega$-continuous functions
	$m_{n}:[0,+\infty)\to[0,+\infty)$ satisfying (\ref{hp:ocont}) with
	the same constant $L$, and such that $m_{n}\to m$ uniformly on
	compact sets.  Let $\{u_{0n},u_{1n}\}\subseteq
	\G_{\varphi,r_{0},\alpha+1/2}(A)\times
	\G_{\varphi,r_{0},\alpha}(A)$ be a sequence converging to
	$(u_{0},u_{1})$ in the same space.
	
	Let finally $T_{1}\in(0,T)$ and $R_{1}>R$ be real numbers such 
	that $R_{1}T_{1}<r_{0}$.
	
	Then we have the following conclusions.
	\begin{enumerate}
		\renewcommand{\labelenumi}{(\arabic{enumi})}
		\item For every $n\in\n$ large enough the Cauchy problem
		(\ref{pbm:h-eq}), (\ref{pbm:h-data}) (with of course $m_{n}$,
		$u_{0n}$, $u_{1n}$ instead of $m$, $u_{0}$, $u_{1}$) has at
		least one solution $u_{n}(t)$ in the space
		\begin{equation}
			C^{1}\left([0,T_{1}];\G_{\varphi,r_{0}-R_{1}t,\alpha}(A)\right)\cap
			C^{0}\left([0,T_{1}];\G_{\varphi,r_{0}-R_{1}t,\alpha+1/2}(A)\right).
			\label{th:dep-reg}
		\end{equation}
	
		\item The sequence $\{u_{n}(t)\}$ is relatively compact in the
		space (\ref{th:dep-reg}).
	
		\item Any limit point of $\{u_{n}(t)\}$ is a solution of the
		limit problem.
	\end{enumerate}
\end{thm}

With minimal technicalities the theory can be extended in order to
allow time-dependent right-hand sides $f_{n}(t)$ with suitable
regularity assumptions, for example in
$L^{2}\left([0,T],\G_{\varphi,r_{0},1/4}(A)\right)$.  We spare the
reader from the details.

\paragraph{Derivative loss}

Our goal is now to prove the optimality of the spaces involved in the 
local existence result. To this end we show that solutions with less 
regular data can exhibit an instantaneous derivative loss. Let us 
introduce the precise notion.

\begin{defn}
	\begin{em}
		Let $H$ and $A$ be as in Theorem~\ref{thm:local}.  Let
		$\varphi:[0,+\infty)\to[1,+\infty)$ be any function, and let
		$\alpha\geq 1/4$.  We say that a solution $u$ of problem
		(\ref{pbm:h-eq}), (\ref{pbm:h-data}) has \emph{instantaneous
		strong derivative loss} of type
		\begin{equation}
			\G_{\varphi,\infty,3/4}(A)\times
			\G_{\varphi,\infty,1/4}(A)\to D(A^{\alpha+1/2})\times
			D(A^{\alpha})
			\label{eq:der-loss}
		\end{equation}
		if the following conditions are fulfilled.
		
		\begin{enumerate}
			\item[(1)] \emph{Regularity of the solution}. There exists 
			$T_{0}>0$ such that 
			\begin{equation}
				u\in C^{1}\left([0,T_{0}];D(A^{1/4})\right) \cap
				C^{0}\left([0,T_{0}];D(A^{3/4})\right).
				\label{hp:u-reg}
			\end{equation}
		
			\item[(2)] \emph{High regularity at $t=0$}. We have that
			$$(u_{0},u_{1})\in
			\G_{\varphi,\infty,3/4}(A)\times\G_{\varphi,\infty,1/4}(A),$$
	
			\item[(3S)] \emph{Low regularity for subsequent times}.
			We have that 
			$$(u(t),u'(t))\not\in
			D(A^{\alpha+1/2+\ep})\times D(A^{\alpha+\ep})
			\quad\quad\forall \ep>0,\ \forall t\in(0,T_{0}].$$ 
		\end{enumerate}
		
		We say that the same solution has \emph{instantaneous weak
		derivative loss} of type~(\ref{eq:der-loss}) if it satisfies
		(1), (2), and
		\begin{enumerate}
			\item[(3W)] \emph{Unboundedness as $t\to 0^{+}$}. There 
			exists a sequence $\tau_{k}\to 0^{+}$ such that
			$$\left|A^{\alpha+1/2+\ep}u(\tau_{k})\right|
			\to +\infty
			\quad\quad\forall\ep>0.$$
		\end{enumerate}
	\end{em}
\end{defn}

The second notion is weaker in the sense that what is actually lost is
the control on the norm of $(u(t),u'(t))$ in
$D(A^{\alpha+1/2+\ep})\times D(A^{\alpha+\ep})$ as $t\to 0^{+}$.  We
are now ready to state our counterexamples, the first one in the
strictly hyperbolic case, the second one in the weakly hyperbolic
case.

\begin{thm}[Derivative loss: strictly hyperbolic case]\label{thm:counterex-ndg}
	Let $A$ be a self-adjoint linear operator on a Hilbert space $H$.
	Let us assume that there exist a countable (not necessarily
	complete) orthonormal system $\{e_{k}\}_{k\geq 1}$ in $H$, and an
	increasing unbounded sequence $\{\lambda_{k}\}_{k\geq 1}$ of
	positive real numbers such that $Ae_{k}=\lambda_{k}^{2}e_{k}$ for
	every $k\geq 1$.
	
	Let $\omega:[0,+\infty)\to[0,+\infty)$ be a continuity modulus 
	such that $\sigma\to\sigma/\omega(\sigma)$ is a nondecreasing 
	function.
	
	Let $\varphi:[0,+\infty)\to[1,+\infty)$ be a function such that
	\begin{equation}
		\lim_{k\to +\infty}\frac{\lambda_{k}}{\varphi(\lambda_{k})}
		\omega\left(\frac{1}{\lambda_{k}}\right)=+\infty.
		\label{hp:derloss-ndg}
	\end{equation}
	
	Then there exist an $\omega$-continuous function
	$m:[0,+\infty)\to[1/2,3/2]$, and a solution $u$ of the
	corresponding problem (\ref{pbm:h-eq}) with
	instantaneous strong derivative loss of type
	$$\G_{\varphi,\infty,3/4}(A)\times \G_{\varphi,\infty,1/4}(A)\to
	D(A^{3/4})\times D(A^{1/4}).$$
\end{thm}

\begin{thm}[Derivative loss: weakly hyperbolic case]\label{thm:counterex-dg}
	Let $H$, $A$, $\{e_{k}\}$, $\{\lambda_{k}\}$, $\omega$ be as in
	Theorem~\ref{thm:counterex-ndg}.  	
	Let $\varphi:[0,+\infty)\to[1,+\infty)$ be a function such that
	\begin{equation}
		\lim_{k\to +\infty}\lambda_{k}
		\left[\varphi\left(\frac{\lambda_{k}}{
		\sqrt{\omega(1/\lambda_{k})}}\right)\right]^{-1}=+\infty.
		\label{hp:derloss-dg}
	\end{equation}
	
	Then there exist an $\omega$-continuous function
	$m:[0,+\infty)\to[0,3/2]$, and a solution $u$ of the
	corresponding equation (\ref{pbm:h-eq}) with instantaneous strong
	derivative loss of type
	$$\G_{\varphi,\infty,3/4}(A)\times \G_{\varphi,\infty,1/4}(A)\to
	D(A)\times D(A^{1/2}),$$
	and instantaneous weak derivative loss of type
	$$\G_{\varphi,\infty,3/4}(A)\times \G_{\varphi,\infty,1/4}(A)\to
	D(A^{3/4})\times D(A^{1/4}).$$
\end{thm}

Note that assumptions (\ref{hp:derloss-ndg}) and (\ref{hp:derloss-dg})
are the counterpart of (\ref{hp:phi-ndg}) and (\ref{hp:phi-dg}),
respectively.  Table~\ref{tab:dl-ndg} and Table~\ref{tab:dl-dg} below
present examples of functions $\omega$ and $\varphi$ satisfying the
assumptions of Theorem~\ref{thm:counterex-ndg} and
Theorem~\ref{thm:counterex-dg} above.  They are the counterpart of
Table~\ref{tab:local-sh} and Table~\ref{tab:local-wh}, respectively.

These examples are based on the construction introduced in \cite{dgcs}
and \cite{cjs} in the linear context.  In those papers the authors
gave examples of coefficients $c(t)$ and initial data $(u_{0},u_{1})$
in such a way that the ``solution'' of the linear problem has
instantaneous derivative loss 
$$\G_{\varphi,\infty,1/2}(A)\times \G_{\varphi,\infty,0}(A)\to
\mbox{hyperdistributions}.$$

This means that the solution is quite regular at time $t=0$, but it is
even outside the space of distributions for $t>0$.  This is usually
presented as a nonexistence result in the space of distributions, but
it is proved by showing that the solution exists and is unique (since
the equation is linear) in a space of hyperdistributions (which in our
notations is a space of the form $\G_{\varphi,r,0}$ with $r<0$), but
exhibits an instantaneous derivative loss up to hyperdistributions.
In the quoted papers the derivative loss is always intended in the
weak sense, but those examples are quite flexible and can be modified
in order to obtain the derivative loss even in the strong sense.

\begin{table}[htb]
	\centering
	\renewcommand{\arraystretch}{1.4}
	\begin{tabular}{|c|c||c|c|}
		\hline
		$m(\sigma)=$ & $m$ is \ldots & $\varphi(\sigma)=$ &
		Derivative loss for data in \ldots \\
		\hline\hline
		\rule[-3ex]{0ex}{7ex}$\displaystyle{\frac{1}{|\log\sigma|^{1/2}}}$ & just continuous &
		$\displaystyle{\frac{\sigma}{\log\sigma}}$ & quasi-analytic
		functions \\
		\hline
		$\sigma^{\beta}$ & $\beta$-H\"{o}lder cont.  & 
		\rule[-3ex]{0ex}{7.5ex}
		$\displaystyle{\frac{\sigma^{1-\beta}}{\log\sigma}}$ & Gevr.\
		sp.\ of order $>(1-\beta)^{-1}$ \\
		\hline
		$\sigma|\log\sigma|^{3}$ & $\beta$-H\"{o}ld.\ cont.\
		$\forall\beta\in(0,1)$ & $\log^{2}\sigma$ & $D(A^{\infty})$ \\
		\hline
	\end{tabular}
	\caption{Pairs of functions $m$, $\varphi$ for which a derivative
	loss example can be found in the strictly hyperbolic case}\medskip
	\label{tab:dl-ndg}
\end{table}

\begin{table}[htb]
	\centering
	\renewcommand{\arraystretch}{1.4}
	\begin{tabular}{|c|c||c|c|}
		\hline
		$m(\sigma)=$ & $m$ is \ldots & $\varphi(\sigma)=$ &
		Derivative loss for data in \ldots \\
		\hline\hline
		$\sigma^{\beta}$ & $\beta$-H\"{o}lder cont.  &
		\rule[-3ex]{0ex}{7.5ex}
		$\displaystyle{\frac{\sigma^{2/(\beta+2)}}{\log\sigma}}$ &
		Gevrey sp.\ of order $>1+\beta/2$ \\
		\hline
		$\sigma$ & Lipschitz cont.  & \rule[-3ex]{0ex}{7.5ex}
		$\displaystyle{\frac{\sigma^{2/3}}{\log\sigma}}$&
		Gevrey sp.\ of order $>3/2$ \\
		\hline
	\end{tabular}
	\caption{Pairs of functions $m$, $\varphi$ for which a derivative
	loss example can be found in the weakly hyperbolic case}\medskip
	\label{tab:dl-dg}
\end{table}

Our strategy is similar.  In the first step we modify the parameters
in those examples in order to stop the derivative loss up to the
$D(A^{3/4})\times D(A^{1/4})$ level.  In the second step we find a
function $m$ in such a way that the coefficient $c(t)$ is actually
equal to $\m{u(t)}$.  This can be easily done as soon as the function
$t\to|A^{1/2}u(t)|^{2}$ is invertible in a neighborhood of $t=0$, and
this can be obtained by modifying just one dominant component of
$u(t)$.  We refer to \cite{gg:local} for the details.

\setcounter{equation}{0}
\section{Uniqueness}\label{sec:uniqueness}

\paragraph{Previous works}

As one can easily guess, uniqueness holds whenever $m$ is (locally)
Lipschitz continuous.  In the strictly hyperbolic case a proof of this
result is contained for example in~\cite{ap}, of course for initial
data in $D(A^{3/4})\times D(A^{1/4})$.  In the weakly hyperbolic case
a proof of the same result is given in~\cite{as} for analytic initial
data.  Now from Theorem~\ref{thm:local} we know that, when $m$ is
Lipschitz continuous and the equation is degenerate, local solutions
exist for all initial data in
$\G_{\varphi,r_{0},3/4}(A)\times\G_{\varphi,r_{0},1/4}(A)$ with
$\varphi(\sigma)=\sigma^{2/3}$.  The uniqueness result under these
assumptions has never been put into writing, but it can be easily
proved by standard arguments.  The main tool is indeed always the
same, namely a Gronwall type lemma for the difference between two
solutions.

As a general fact, uniqueness for a nonlinear evolution equation is
much more difficult to establish if the nonlinear term is not locally
Lipschitz continuous.  Therefore it is hardly surprising that also in
the case of Kirchhoff equations the non-Lipschitz case remained widely
unexplored for a long time.  To our knowledge indeed uniqueness issues
have been previously considered only in section~4 of~\cite{as}, where
two results are presented.

The first one is a one-dimensional example ($H=\re$) where problem
(\ref{pbm:h-eq}), (\ref{pbm:h-data}) admits infinitely many local
solutions.  The second result is a detailed study of the case where
$u_{0}$ and $u_{1}$ are \emph{eigenvectors} of $A$ relative to the
\emph{same eigenvalue}.  In this special situation (which can be
easily reduced to the two dimensional case $H=\re^{2}$) the authors proved
that uniqueness of the local solution fails if and only if the
following three conditions are satisfied:
\begin{list}{}{\leftmargin 4em \labelwidth 4em}
	\item[(AS1)] $\langle Au_{0},u_{1}\rangle= 0$,

	\item[(AS2)] $|A^{1/2}u_{1}|^{2}-\m{u_{0}}|Au_{0}|^{2}=0$,

	\item[(AS3)]  $m$ satisfies a suitable integrability condition in a
	neighborhood of $|A^{1/2}u_{0}|^{2}$.
\end{list}

As a consequence, the local solution is unique if at least one of the
conditions above is not satisfied.

\paragraph{Our uniqueness result}

Our contribution is the extension of the first two parts of the above
result from the two dimensional case with equal eigenvalues to the
infinite dimensional case with arbitrary eigenvalues.  In other words,
we prove that in the general case the solution is necessarily unique
whenever either (AS1) or (AS2) are not satisfied.  The precise
statement is the following.

\begin{thm}[Uniqueness]\label{thm:uniqueness}
	Let $H$, $A$, $\omega$, $m$, $\varphi$, $\Lambda$ be as in
	Theorem~\ref{thm:local}.  Let us assume that
	\begin{equation}
		(u_{0},u_{1})\in\G_{\varphi,r_{0},3/2}(A)\times
		\G_{\varphi,r_{0},1}(A) 
		\label{hp:main-data}
	\end{equation}
	for some $r_{0}>0$, and
	\begin{equation}
		\left|\langle Au_{0},u_{1}\rangle\right|+
		\left||A^{1/2}u_{1}|^{2}-\m{u_{0}}|Au_{0}|^{2}\right|\neq 0.
		\label{hp:main}
	\end{equation}
	
	Let us assume that problem (\ref{pbm:h-eq}), (\ref{pbm:h-data})
	admits two local solutions $v_{1}$ and $v_{2}$ in
	\begin{equation}
		C^{2}\left([0,T];\G_{\varphi,r_{1},1/2}(A)\right)\cap
		C^{1}\left([0,T];\G_{\varphi,r_{1},1}(A)\right)\cap
		C^{0}\left([0,T];\G_{\varphi,r_{1},3/2}(A)\right)
		\label{hp:reg-sol}
	\end{equation}
	for some $T>0$, and some $r_{1}\in(0,r_{0})$.
	
	Then we have the following conclusions.
	\begin{enumerate}
		\renewcommand{\labelenumi}{(\arabic{enumi})}
		\item  There exists $T_{1}\in(0,T]$ such that
		\begin{equation}
			v_{1}(t)=v_{2}(t)
			\hspace{2em}
			\forall t\in[0,T_{1}].
			\label{th:uniqueness}
		\end{equation}
	
		\item  Let $T_{*}$ denote the supremum of all $T_{1}\in (0,T]$
		for which (\ref{th:uniqueness}) holds true. Let $v(t)$ denote 
		the common value of $v_{1}$ and $v_{2}$ in $[0,T_{*}]$.
		
		Then either $T_{*}=T$ or
		$$\left|\langle Av(T_{*}),v'(T_{*})\rangle\right|+
		\left||A^{1/2}v'(T_{*})|^{2}-
		\m{v(T_{*})}|Av(T_{*})|^{2}\right|= 0.$$
	\end{enumerate}
\end{thm}

Let us make some comments on the assumptions.  Inequality
(\ref{hp:main}) is equivalent to say that either (AS1) or (AS2) are
not satisfied.  The space (\ref{hp:reg-sol}) is the natural one when
initial data satisfy (\ref{hp:main-data}).  Indeed from the
propagation of regularity (see Theorem~\ref{thm:apriori-est}) it
follows that any solution $u(t)$ satisfying (\ref{hp:reg-u}) with
initial data as in (\ref{hp:main-data}) lies actually in
(\ref{hp:reg-sol}).  Assumption (\ref{hp:main-data}) on the initial
data is stronger than the corresponding assumption in
Theorem~\ref{thm:local}.  This is due to a technical point in the
proof.  However in most cases the difference is only apparent.  For
example if $\omega(\sigma)=\sigma^{\beta}$ for some $\beta\in(0,1]$,
then the following implication $$u\in\G_{\varphi,r,0}(A)\
\Longrightarrow\
u\in\G_{\varphi,r-\ep,\alpha}(A)$$
holds true for every $r>0$, $\ep\in(0,r)$, $\alpha\geq 0$.
Therefore in this case every solution satisfying
(\ref{th:reg-sol}) fulfils (\ref{hp:reg-sol}) with
$r_{1}=(r_{0}-RT)/2$.

In the proof of Theorem~\ref{thm:uniqueness}, for which we refer
to~\cite{gg:uniqueness}, we introduced a technique which seems to be
new, and hopefully useful to handle also different evolution equations
with non-Lipschitz terms.  The main idea is to split the uniqueness
problem in two steps, which we call trajectory uniqueness and
parametrization uniqueness.

\subparagraph{\emph{\textmd{Trajectory uniqueness}}}

The first step of the proof consists in showing that the image of the
curve $(A^{1/2}u(t),u'(t))$ in the phase space (for example in
$D(A^{3/4})\times D(A^{1/4})$) is unique.  To this end we introduce
the new variable
$$s=\psi(t):=|A^{1/2}u(t)|^{2}-|A^{1/2}u_{0}|^{2}.$$
	
If the function $\psi$ is invertible in a right-hand neighborhood of 
the origin, then we can parametrize the curve using the variable $s$. 
If $(z(s),w(s))$ is this new parametrization, then $z$ and $w$ are 
solutions of the following system
\begin{equation}
	z'(s)=\frac{A^{1/2}w(s)}{2\langle A^{1/2}z(s),w(s)\rangle},
	\quad\quad
	w'(s)=-m\left(s+|A^{1/2}u_{0}|^{2}\right) \frac{A^{1/2}z(s)}{2\langle
	A^{1/2}z(s),w(s)\rangle},
	\label{system}
\end{equation}
with initial data
\begin{equation}
	z(0)=A^{1/2}u_{0},
	\quad\quad
	w(0)=u_{1}.
	\label{system:data}
\end{equation}

What is important is that the non-Lipschitz term $m(|A^{1/2}u|^{2})$
of the original equation has become the non-Lipschitz coefficient
$m(s+|A^{1/2}u_{0}|^{2})$ in the second equation of system
(\ref{system}), and it is well known that nonregular
\emph{coefficients} do not affect uniqueness.  Therefore the solution
of the system is unique.

\subparagraph{\emph{\textmd{Parametrization uniqueness}}}

The second part of the proof consists in showing that the unique 
trajectory obtained in the previous step can be covered by solutions 
in a unique way. To this end, we first show that the parametrization 
$\psi(t)$ is a solution of the Cauchy problem
\begin{equation}
	\psi'(t)=F(\psi(t)),
	\quad\quad
	\psi(0)=0,
	\label{pbm:cauchy}
\end{equation}
where $F(\sigma):=2\langle A^{1/2}z(\sigma),w(\sigma)\rangle$.  The
function $F$ is just continuous in $\sigma=0$, and this in not enough
to conclude that the solution of (\ref{pbm:cauchy}) is unique.  On the
other hand, the differential equation in (\ref{pbm:cauchy}) is
\emph{autonomous}, and for autonomous equations it is well known that
there is a unique solution such that $\psi(t)>0$ for $t>0$.

Proving that $\psi(t)>0$ for $t>0$, and more generally that $\psi$ is
invertible in a right-hand neighborhood of $t=0$ (as required in the
first step), is the point where the quite strange assumptions (AS1)
and (AS2) play their role.  Indeed we have that
\begin{eqnarray*}
	\psi'(0)=0 & \Longleftrightarrow & \mbox{(AS1) holds true},  \\
	\psi''(0)=0 & \Longleftrightarrow & \mbox{(AS2) holds true}.
\end{eqnarray*}

If (\ref{hp:main}) is true, then either (AS1) or (AS2) are false, 
hence either $\psi'(0)\neq 0$ or $\psi''(0)\neq 0$. In both cases 
$\psi(t)$ is invertible where needed.

We conclude by pointing out that the denominators in (\ref{system})
are actually $\psi'(t)$, hence they can vanish for $t=0$.  Since in
that case we have that $\psi''(0)\neq 0$, then for sure denominators
are different from $0$ for all $t>0$ small enough, and their vanishing
in $t=0$ is of order one.  This kind of singularity doesn't affect
existence or uniqueness for system (\ref{system}), but it is in some
sense the limit exponent.  For this reason we cannot deal with the
same technique the case where $\psi'(0)=\psi''(0)=0$ but
$\psi'''(0)\neq 0$, which originates denominators with a singularity
of order 2.

\setcounter{equation}{0}
\section{Global existence}\label{sec:global}

\paragraph{Previous works}

Global existence for Kirchhoff equations has been proved in at least 
five special cases.

\subparagraph{\emph{\textmd{Analytic data}}}

This is the result we quoted as (B) in the history of local existence 
results. We recall the main assumptions: the equation is weakly 
hyperbolic, the nonlinearity is continuous, and initial data are 
analytic.

\subparagraph{\emph{\textmd{Quasi-analytic data}}}

K.\ Nishihara~\cite{nishihara} proved global existence for
a class of initial data which strictly contains analytic functions.
His assumptions are that the equation is
strictly hyperbolic, the nonlinearity is Lipschitz continuous, and
$$(u_{0},u_{1})\in\G_{\varphi,r_{0},1/2}(A)\times
\G_{\varphi,r_{0},0}(A),$$
where $r_{0}>0$, and $\varphi:[0,+\infty)\to[1,+\infty)$ is an
increasing function satisfying suitable convexity and integrability
conditions.  He proves existence of a global solution  
$$u\in
C^{1}\left([0,+\infty);\G_{\varphi,r_{0},0}(A)\right)\cap
C^{0}\left([0,+\infty);\G_{\varphi,r_{0},1/2}(A)\right).$$

We point out that, in contrast with our local existence results, this
solution lives in a Hilbert \emph{space}, instead of a Hilbert
\emph{scale}.

The most celebrated example of function $\varphi$ satisfying the
assumptions is $\varphi(\sigma)=\sigma/\log\sigma$, in which case one
has global existence in a space which contains non-analytic initial
data.  On the contrary, the function $\varphi(\sigma)=\sigma^{\beta}$
with $\beta<1$ never satisfies the assumptions.  In other words,
Nishihara's spaces are intermediate classes between Gevrey and
analytic functions.

It would be interesting to compare Nishihara's assumptions with
$$\int_{1}^{+\infty}\frac{\varphi(\sigma)}{\sigma^{2}}\,d\sigma 
=+\infty,$$
which is the usual definition of quasi-analytic classes.

\subparagraph{\emph{\textmd{Special nonlinearities}}}

In a completely different direction, S.\ I.\
Pohozaev~\cite{poho-m} considered the special case where
$m(\sigma):=(a+b\sigma)^{-2}$ for some $a>0$ and $b\in\re$.  He proved
global existence for initial data $(u_{0},u_{1})\in D(A)\times
D(A^{1/2})$ satisfying the nondegeneracy condition
$a+b|A^{1/2}u_{0}|^{2}>0$.

The main point is that in this case (and in a certain sense only in
this case) equation (\ref{pbm:h-eq}) admits the second order nonnegative
invariant 
$$\mathcal{P}(t):=
\left(a+b|A^{1/2}u(t)|^{2}\right)|A^{1/2}u'(t)|^{2}+
\frac{|Au(t)|^{2}}{a+b|A^{1/2}u(t)|^{2}}- \frac{b}{4}\langle
Au(t),u'(t)\rangle^{2}.$$

Exploiting that $\mathcal{P}(t)$ is constant, it is not difficult to
obtain a uniform bound on $\langle Au(t),u'(t)\rangle^{2}$, from which
global existence follows in a standard way.

Recently, some new results have been obtained along this path.  The
interested reader is referred to \cite{by}.

\subparagraph{\emph{\textmd{Dispersive equations}}}

Global existence results have been obtained for the concrete equation
(\ref{eq:k}) in cases where \emph{dispersion} plays a crucial role,
namely when $\Omega=\re$ (see J.\ M.\ Greenberg and S.\ C.\
Hu~\cite{gh}), $\Omega=\re^{n}$ (see P.\ D'Ancona and S.\
Spagnolo~\cite{das}), or $\Omega=$ exterior domain (see T.\
Yamazaki~\cite{yamazaki1,yamazaki2} and the references quoted
therein).

The prototype of these results is global existence provided that the
equation is strictly hyperbolic, the nonlinearity is Lipschitz
continuous, and initial data have Sobolev regularity and satisfy
suitable \emph{smallness assumptions} and \emph{decay conditions} at
infinity. We refer to the quoted literature for precise statements.

\subparagraph{\emph{\textmd{Spectral gap initial data}}}

More recently, R.\ Manfrin~\cite{manfrin1} (see also
\cite{manfrin2}, \cite{hirosawa2}) proved global existence in a new
class of nonregular initial data. In order to describe the most 
astonishing aspect of his work, we need the following definition.

\begin{defn}
	\begin{em}
		Let $\mathcal{M}$ and $\mathcal{F}$ be two \emph{subsets} of
		$D(A^{1/2})\times H$.  We say that $\mathcal{M}$ has the ``Sum
		Property'' in $\mathcal{F}$ if
		$\mathcal{M}\subseteq\mathcal{F}$ and
		$\mathcal{M}+\mathcal{M}\supseteq\mathcal{F}$.

		In other words, for every $(u_{0},u_{1})\in\mathcal{F}$ there exist
		$(\overline{u}_{0},\overline{u}_{1})\in\mathcal{M}$ and
		$(\widehat{u}_{0},\widehat{u}_{1})\in\mathcal{M}$ such that
		$u_{0}=\overline{u}_{0}+\widehat{u}_{0}$, and
		$u_{1}=\overline{u}_{1}+\widehat{u}_{1}$.
	\end{em}
\end{defn}

Let us assume now that the equation is strictly hyperbolic, and the
nonlinearity is of class $C^{2}$.  The main result of~\cite{manfrin1}
is that there exists a subset $\mathcal{M}\subseteq D(A)\times
D(A^{1/2})$ such that $\mathcal{M}$ has the ``Sum Property'' in
$D(A)\times D(A^{1/2})$, and problem (\ref{pbm:h-eq}),
(\ref{pbm:h-data}) admits a global solution for every
$(u_{0},u_{1})\in\mathcal{M}$.

As a corollary, any initial condition $(u_{0},u_{1})\in D(A)\times
D(A^{1/2})$ is the sum of two pairs of initial conditions for which
the solution is global! Of course the set $\mathcal{M}$ is not a 
vector space, but just a star-shaped subset (actually a cone). We 
refer to the quoted papers for the definition of $\mathcal{M}$.

\paragraph{Our global existence result}

In \cite{gg:global} we proved a result in the same spirit of Manfrin's
one, but without assuming the strict hyperbolicity or the regularity
of the nonlinearity.

Let $\mathcal{L}$ denote the set of all sequences $\{\rho_{n}\}$ of
positive real numbers such that $\rho_{n}\to +\infty$.  Let
$\varphi:[0,+\infty)\to[1,+\infty)$ be any function, and let
$\alpha\geq 0$, $\beta\geq 0$. We can now define what we call 
\emph{generalized Gevrey-Manfrin spaces}, namely
\begin{equation}
	\GM_{\varphi,\{\rho_{n}\},\alpha}^{(\beta)}(A):=\left\{ u\in
	H:\sum_{\lk>\rho_{n}}\lambda_{k}^{4\alpha}u_{k}^{2}
	\exp\left(\rho_{n}^{\beta}\varphi(\lk)\right)\leq\rho_{n} \quad
	\forall n\in\n\right\},
	\label{defn:GM}
\end{equation}
and		
\begin{equation}
	\GM_{\varphi,\alpha}^{(\beta)}(A):=
	\bigcup_{\{\rho_{n}\}\in\mathcal{L}}
	\GM_{\varphi,\{\rho_{n}\},\alpha}^{(\beta)}(A).
	\label{defn:GM-cone}
\end{equation}

Admittedly this definition has no immediate interpretation.  Let us
compare (\ref{defn:GM}) with (\ref{defn:G}) and (\ref{defn:trebar}).
In the inequalities in (\ref{defn:GM}) the weight $\rho_{n}$ appears in
the right-hand side, and in the left-hand side in place of $r$.
Moreover $\rho_{n}$ appears also in the summation, which is now
restricted to eigenvalues $\lambda_{k}>\rho_{n}$.  The weight in
the left-hand side is inside an exponential term, hence it dominates
on the weight in the right-hand side.  It follows that the
inequalities in (\ref{defn:GM}) are smallness assumptions on the
``tails'' of suitable series.  More important, the smallness is not
required for all tails, but only for a subsequence.

It is easy to see that the space defined in (\ref{defn:GM}) is
actually a vector space, while the space defined by
(\ref{defn:GM-cone}) is a cone in $\G_{\varphi,\infty,\alpha}(A)$,
because its elements may be defined starting from different sequences
in $\mathcal{L}$.  This fact is crucial in the proof of the 
``Sum Property'' (see Proposition~3.2 in~\cite{gg:global}).

\begin{prop}[Sum Property]\label{prop:sp}
	For every $\varphi:[0,+\infty)\to[1,+\infty)$, $\alpha\geq 0$,
	$\beta> 0$ we have that 
	$$\GM_{\varphi,\alpha+1/2}^{(\beta)}(A)\times
	\GM_{\varphi,\alpha}^{(\beta)}(A)$$
	has the ``Sum Property'' in 
	$$\G_{\varphi,\infty,\alpha+1/2}(A)\times
	\G_{\varphi,\infty,\alpha}(A).$$
\end{prop}
	
The proof of the ``Sum Property'' is based on the following idea.  Let
us consider an increasing and divergent sequence $s_{n}$ of positive
real numbers.  Then any $u_{0}\in H$ can be written as the sum of
$\overline{u}_{0}$ and $\widehat{u}_{0}$, where $\overline{u}_{0}$ has the same
components of $u_{0}$ with respect to eigenvectors corresponding to
eigenvalues belonging to intervals of the form
$[s_{2n},s_{2n+1})$, and components equal to 0 with respect to the
remaining eigenvectors, and vice versa for $\widehat{u}_{0}$.  If the
sequence $s_{n}$ grows fast enough, then it turns out that
$\overline{u}_{0}$ and $\widehat{u}_{0}$ lie in suitable generalized
Gevrey-Manfrin spaces corresponding to the sequences $s_{2n}$ and
$s_{2n+1}$.  Note that the spectrum of both $\overline{u}_{0}$ and
$\widehat{u}_{0}$ has a sequence of ``big holes'', which justify the term
``spectral gap'' initial data.

We are now ready to state our global existence result (see Theorem~3.1
and Theorem~3.2 in~\cite{gg:global}).

\begin{thm}[Global existence]\label{thm:global}
	Let $H$, $A$, $\omega$, $m$, $\varphi$, $\Lambda$ be as in
	Theorem~\ref{thm:local}. Let $\{\rho_{n}\}\in\mathcal{L}$, and let
	$$(u_{0},u_{1})\in
	\GM_{\varphi,\{\rho_{n}\},3/4}^{(\beta)}(A)\times
	\GM_{\varphi,\{\rho_{n}\},1/4}^{(\beta)}(A),$$
	where $\beta=2$ if the equation is strictly hyperbolic, and 
	$\beta=3$ if the equation is weakly hyperbolic.
			
	Then problem (\ref{pbm:h-eq}), (\ref{pbm:h-data}) admits at least
	one global solution $u$ with 
	$$u\in
	C^{1}\left([0,+\infty);\G_{\varphi,r,3/4}(A)\right)\cap
	C^{0}\left([0,+\infty);\G_{\varphi,r,1/4}(A)\right)$$
	for every $r>0$.
\end{thm}

Combining Theorem~\ref{thm:global} and Proposition~\ref{prop:sp} we
obtain the following statement: every pair of initial conditions
satisfying (\ref{hp:data}) with $r_{0}=\infty$ is the sum of two
pairs of initial conditions for which the solution is global.  We have
thus extended to the general case the astonishing aspect of Manfrin's
result.

The extra requirement that $r_{0}=\infty$ is hardly surprising.  It is
indeed a necessary condition for existence of global solutions even in
the theory of linear equations with nonsmooth time-dependent
coefficients.

We conclude by remarking that, in the concrete case, these spaces do
not contain any compactly supported function.

\setcounter{equation}{0}
\section{Open problems}\label{sec:open}

The main open problem in the theory of Kirchhoff equations is for 
sure the existence of global solutions in $C^{\infty}$. In the 
abstract setting it can be stated as follows.

\begin{open}
	Let us assume that equation (\ref{pbm:h-eq}) is strictly
	hyperbolic, that $m\in C^{\infty}(\re)$, and $(u_{0},u_{1})\in
	D(A^{\infty})\times D(A^{\infty})$.
	
	Does problem (\ref{pbm:h-eq}), (\ref{pbm:h-data}) admit a global 
	solution?
\end{open}

The same problem can be restated in all situations where a local
solution has been proved to exist (see Theorem~\ref{thm:local}).  Up
to now indeed we know no example of local solution, with any
regularity, which is not global.

Now we would like to mention some other open questions.  The first one
concerns local (but of course also global) existence for initial data
in $D(A^{1/2})\times H$, which is the natural \emph{energy space} for
a second order wave equation.

\begin{open}
	Let us assume that equation (\ref{pbm:h-eq})
	is strictly hyperbolic, that $m\in C^{\infty}(\re)$, and
	$(u_{0},u_{1})\in D(A^{\alpha+1/2})\times D(A^{\alpha})$ for some 
	$\alpha\in[0,1/4)$.
	
	Does problem (\ref{pbm:h-eq}), (\ref{pbm:h-data}) admit a local 
	solution? Of course in this case we accept solutions 
	\begin{equation}
		u\in C^{1}([0,T_{0}];H) \cap C^{0}([0,T_{0}];D(A^{1/2})).
		\label{reg:es}
	\end{equation}
\end{open}

Once again we know no counterexample, even with degenerate equations 
or nonlinearities which are just continuous.

We stress that counterexamples are the missing element in all the
theory.  We proved the optimality of our local existence results by
showing examples of solutions with derivative loss.  These are
actually counterexamples to propagation of regularity, but not
counterexamples to existence.  We can therefore ask the following
question.

\begin{open}
	Do there exist a nonnegative continuous function $m$, and initial
	data $(u_{0},u_{1})\in D(A^{3/4})\times D(A^{1/4})$, such that
	problem (\ref{pbm:h-eq}), (\ref{pbm:h-data}) admits no (local)
	solution $u$ satisfying (\ref{hp:u-reg})?
	
	Do there exist a nonnegative continuous function $m$, and initial
	data $(u_{0},u_{1})\in D(A^{1/2})\times H$, such that
	problem (\ref{pbm:h-eq}), (\ref{pbm:h-data}) admits no (local)
	solution $u$ satisfying (\ref{reg:es})?
\end{open}

We conclude by mentioning three open questions related to uniqueness
results.  The first one concerns once again counterexamples.  The
motivation is that we know no example where uniqueness fails apart
from those given in~\cite{as}.  So we ask whether different
counterexamples can be provided.

\begin{open}
	Let $H$, $A$, $\omega$, $m$, $\varphi$, $A$, $u_{0}$, $u_{1}$ be
	as in Theorem~\ref{thm:uniqueness}, but without assumption
	(\ref{hp:main}).  Let us assume that problem (\ref{pbm:h-eq}),
	(\ref{pbm:h-data}) admits two local solutions.
	
	Can one conclude that $u_{0}$ and $u_{1}$ are eigenvectors of $A$
	relative to the same eigenvalue?
\end{open}

We stress that this problem is open even in the simple case
$H=\re^{2}$, where $\omega$ and $\varphi$ play non role, and no
regularity is required on initial data.

The second open problem concerns trajectory uniqueness, namely the key
step in the proof of our uniqueness result.  One can indeed observe
that, even in the non-uniqueness examples of \cite{as}, all the
different solutions describe (a subset of) the same trajectory with a
different pace.  With our notations this is equivalent to say that the
solution of (\ref{system}), (\ref{system:data}) is unique.  We ask
whether this property is true in general.

\begin{open}
	Let $H$, $A$, $\omega$, $m$, $\varphi$, $A$, $u_{0}$, $u_{1}$ be
	as in Theorem~\ref{thm:uniqueness}, but without assumption
	(\ref{hp:main}).  Let us consider system (\ref{system}), with
	initial data (\ref{system:data}).
	
	Does this system admit at most one solution?
\end{open}

Note that in the case where $\langle Au_{0},u_{1}\rangle=0$ it is by
no means clear that the system admits at least one solution, since
this implicitly requires that $\langle A^{1/2}z(s),w(s)\rangle\neq 0$
for every $s\in(0,s_{0}]$.  In any case the above question doesn't
concern existence, but just uniqueness provided that a solution
exists.

The last open problem concerns the regularity assumptions on initial
data and solutions required in the uniqueness result.  Indeed in
Theorem~\ref{thm:uniqueness} we proved that inequality (\ref{hp:main})
yields uniqueness provided that initial data satisfy
(\ref{hp:main-data}) and solutions satisfy (\ref{hp:reg-sol}).
Similar assumptions are required in the uniqueness result for
Lipschitz continuous nonlinearities.  On the other hand, solutions of
problem (\ref{pbm:h-eq}), (\ref{pbm:h-data}) may exist also if
(\ref{hp:main-data}) is not satisfied (this is the case, for example,
of our solutions with derivative loss).  We ask whether uniqueness
results can be proved for these solutions.

\begin{open}
	Is it possible to prove the known uniqueness results (namely the
	Lipschitz case and our Theorem~\ref{thm:uniqueness}) with less
	regularity requirements on initial data or for solution in the
	energy space?
\end{open}

Just to give an extreme example, let us consider problem
(\ref{pbm:h-eq}), (\ref{pbm:h-data}) in the strictly hyperbolic case,
with an analytic nonlinearity $m$, and analytic initial data.  We know
that there exists a unique solution in $D(A^{3/4})\times D(A^{1/4})$,
which is actually analytic.  However, as far as we know, no one can
exclude that there exists a different solution in $D(A^{1/2})\times H$
with the same initial data!

\subsubsection*{\centering Acknowledgments}

This note is an extended version of the talk presented by the second
author in the section ``Dispersive Equations'' of the 7th ISAAC
conference (London 2009).  We would like to thank once again the
organizers of that section, Prof.\ F.\ Hirosawa and Prof.\ M.\
Reissig, for their kind invitation.

\label{NumeroPagine}

\end{document}